\documentclass{article}
\usepackage{amssymb}
\usepackage{amsmath}
\title{ Vector ultrametric spaces and a fixed point theorem
for correspondences \\[0.3cm]}

\author{{Kourosh Nourouzi  \thanks {e-mail: nourouzi@kntu.ac.ir; fax: +98 21
22853650}
}\\[0.4cm]
{  Department of Mathematics,  K. N. Toosi University of Technology,}\\
{\em P.O. Box 16315-1618, Tehran, Iran.}\\
}

\newtheorem{definition}{Definition}
\newtheorem{corollary}{Corollary}
\newtheorem{ex}{Example}

\begin{document}

\maketitle \begin{abstract} In this paper,  vector ultrametric
spaces are introduced and   a fixed point  theorem is given for
correspondences.  Our main result generalizes a known theorem in
ordinary ultrametric spaces.

\end{abstract}

\renewcommand{\baselinestretch}{1.1}
\def\thefootnote{ \ }

\footnotetext{{\em} $2010$ Mathematics Subject Classification.
Primary: 26E30, Secondary: 47H10
}

\section{Introduction and Preliminaries}

An \textit{ultrametric space} $(X,d)$ is a metric space in which
the triangle inequality is replaced by
 $$d(x,y)\leq \max \{d(x,z),d(z,y)\}, \,\,\,\,\,\,\,\,\,\ (x,y,z\in X).$$
 A  generalization of  the notion of  ultrametric space via
partially ordered sets was given in  \cite{rib, rib1} which  led
 some applications  to logic programming \cite{rib2}, computational
logic \cite{seda}, and quantitative domain theory \cite{domain}.

    In this paper  we allow ultrametrics to take values in an arbitrary
cone of a complete modular space.  The main result of this paper is
 a fixed point theorem for correspondences  in vector ultrametric spaces
which  generalizes the main theorem  presented in
\cite{petalas}.

We first present some basic notions which will be needed in this paper.

A modular on a real linear space $\mathcal{A}$ is a real valued functional  $\rho$ on  $\mathcal{A}$ which satisfies the  conditions:\\
\indent 1. $\rho(x)=0$ if and only if $x=0$,\\
\indent 2. $\rho(x)=\rho(-x)$,\\
\indent 3. $\rho(\alpha x+\beta y)\leq \rho(x)+\rho(y)$, for all
$x,y\in
\mathcal{A}$ and $\alpha, \beta \geq 0$, $\alpha +\beta=1$. \\

\noindent Then, the  vector subspace $\mathcal{A}_\rho = \{x\in X:
\rho(\alpha x)\rightarrow 0 \,\,\,\rm{as}
\,\,\,\alpha\rightarrow0\}$ of $\mathcal{A}$ is called a modular
space.

The modular $\rho$ is called convex (see, e.g., \cite{agh, ka} for a more general form of convexity) if Condition (3) is replaced with
\begin{center}  $\rho(ax+by)\leq a\rho(x)+b\rho(y)$ for all $x,y\in X$ and all $a,b\geq0$ with $a+b=1$.
\end{center}

\noindent A sequence $(x_n)_{n=1}^\infty$ in $\mathcal{A}_{\rho}$
is called $\rho$-convergent (briefly, convergent) to $x\in
\mathcal{A}_{\rho}$  if $\rho(x_n-x)\rightarrow 0$ as
$n\rightarrow\infty$; $(x_n)_{n=1}^\infty$ is said to be a Cauchy
sequence if $\rho(x_m-x_n)\rightarrow 0$ as
$m,n\rightarrow\infty$. By a $\rho$-closed (briefly, closed) set
in $\mathcal{A}_\rho$ it is meant that it contains the limit of
all its convergent sequences. And, $\mathcal{A_\rho}$ is a
complete modular space if every Cauchy sequence in
$\mathcal{A}_\rho$ is convergent to a point of
$\mathcal{A}_\rho$. The modular $\rho$ is said to satisfy the
$\Delta_2$-condition if there exists  $k>0$ such that
$\rho(2x)\leq k\rho(x)$ for all $x\in \mathcal{A}_\rho$. The reader is referred  to  \cite{koz, mus} for more details, We also suggest the reader see \cite{k1, k2, k5, kl, k3,  k4}.

%


\begin{definition}{\rm A nonempty subset
$\mathcal{P}$ of a  complete modular space  $\mathcal{A}_\rho$ is
called a \textit{cone} if
\begin{description}
\item (i) $\mathcal{P}$ is $\rho$-closed, and
$\mathcal{P}\neq\{0\}$; \item(ii) $a,b\in\mathbb R$ ,$
\;a,b\geq0$,$ \;x,y\in \mathcal{P}\Rightarrow ax+by\in
\mathcal{P}$; \item(iii) $ \mathcal{P}\cap (-\mathcal{P})=\{0\}.$
\end{description}}
\end{definition}

A partial order $\preceq $ can be induced on $\mathcal{A}_\rho$ by
every cone $\mathcal{P}\subset \mathcal{A}$ as  $x\preceq  y$
whenever $y-x\in \mathcal{P}$. A  cone $\mathcal{P}$ is
called \textit{normal} (or $\rho$-normal) if there is a
positive real number $c$ (normal constant) such that
$$0\preceq   x\preceq   y \Rightarrow \,\,\, \rho(x)\leq  c\rho(y),\,\,\,\,\,\,\,\,\,\,\,\,\, (x,y \in \mathcal{A}_\rho).$$

%
%
%
When the modular $\rho$ of  $\mathcal{A}_\rho$ satisfies
$\Delta_2$-condition with $\Delta_2$-constant $k$, it can be
replaced with an equivalent  modular
$\sigma$ satisfying $\Delta_2$-condition  for which the normal
constant of $\mathcal{P}$ is $1$ with respect to $\sigma$. In
fact, for such modular $\rho$ it suffices to define
$$\sigma(x)=\inf_{y\preceq x}\rho(y) + \inf_{x\preceq z}\rho(z)
\,\,\,\,\,\,\,\,\,\,\ (x\in \mathcal{A}_\rho).$$ Then, $\sigma$
is a modular on $\mathcal{A}_\rho$ which is equivalent to $\rho$
and satisfies $\Delta_2$-condition. To see this, we just show
that $x=0$ if $\rho(x)=0$ and $\rho(\alpha x+\beta y)\leq
\rho(x)+\rho(y)$ as $\alpha, \beta \geq 0$, $\alpha +\beta=1$.
Let $\varepsilon>0$ be given. There exist $y,z\in
\mathcal{A}_\rho$  such that $y\preceq x \preceq z$ and
$\max\{\rho(y), \rho(z)\}\leq \varepsilon$. Since $x-y\preceq
z-y$, we get $$\rho(\frac{x}{4})\leq
\rho(\frac{x-y}{2})+\rho(\frac{y}{2})\leq
c\rho(\frac{z-y}{2})+\rho(\frac{y}{2})\leq
c\rho(z)+(c+1)\rho(y),$$ where $c$ is the normal constant. This
implies that  $x=0$. Now let $x,u\in \mathcal{A}_\rho$. Choose
$y_1,y_2,z_1,z_2\in \mathcal{A}_\rho$ such that $y_1\preceq
x\preceq z_1$ and $y_2 \preceq u \preceq z_2$ with
$$\rho(y_1)+\rho(z_1)\leq \sigma(x)+\varepsilon, \,\,\,\,\,\,\
\rho(y_2)+\rho(z_2)\leq \sigma(u)+\varepsilon.$$ Since $\alpha
y_1+\beta y_2\preceq \alpha x+\beta u\preceq \alpha z_1+\beta
z_2$, we have $$\sigma(\alpha x+\beta u)\leq \rho(\alpha
y_1+\beta y_2)+\rho(\alpha z_1+\beta z_2),$$ and consequently
$$\sigma(\alpha x+\beta u)\leq \sigma(x)+\sigma(u)+ 2\varepsilon
.$$ To see  the normal constant of $\sigma$, let $0\preceq x
\preceq u$. Then,  $$\sigma(x)= \inf_{x\preceq z}\rho(z)\leq
\inf_{u\preceq z} \rho(z)=\sigma(u),$$ that is the desired
constant is $1$.  Finally, $\sigma(x)\leq 2 \rho(x)$, for each
$x\in \mathcal{A}_\rho$. On the other hand, if $y\preceq x
\preceq z$, we have $$\rho(\frac{x}{2})\leq
\rho(\frac{x-y}{2})+\rho(\frac{y}{2})\leq
c\rho(\frac{z-y}{2})+\rho(\frac{y}{2}) \leq
(c+1)(\rho(y)+\rho(z)),$$ therefore, $$\rho(\frac{x}{2})\leq
(c+1)\sigma(x).$$ Since $\sigma$ satisfies $\Delta_2$-condition,
we get $$\rho(x)\leq k(c+1)\sigma(x), \,\,\,\,\,\,\,\,\,\,\,\
(x\in \mathcal{A}_\rho).$$

Hence, by a normal cone we always assume that its normal constant
is $1$. We also would say that the cone $P$ is \textit{unital} if
there exists a vector $e\in \mathcal{P}$ with modular $1$ such
that
$$x\preceq \rho(x) e \,\,\,\,\,\,\,\,\,\,\,\ (x\in \mathcal{P}).$$

Throughout this note, we  suppose that $\mathcal{P}$ is a cone in
complete modular space $\mathcal{A}_\rho$ where its modular
is convex and satisfies $\Delta_2$-condition and $\preceq $ is the partial order
induced by
 $\mathcal{P}$.

\begin{definition} {\rm
Let $\mathcal{X}$ be a nonempty set. If the mapping
\mbox{$d:\mathcal{X}\times \mathcal{X}\rightarrow
\mathcal{A_\rho}$} satisfies the following conditions:
\begin{description}
\item(CUM1) $d(x,y)\succeq 0$ for all $x,y\in \mathcal{X}$ and
$d(x,y)=0$ if and only if $x=y$; \item(CUM2) $d(x,y)=d(y,x)$ for
all $x,y\in \mathcal{X}$; \item(CUM3) If $d(x, z)\preceq p$
 and $d(y, z)\preceq p$,
 then $d(x,y)\preceq p$,  for any $x,y,z\in \mathcal{X}$, and $p\in \mathcal{P}$;
\end{description}
then $d$ is called a \textit{vector ultrametric} on $\mathcal{X}$,
and the triple $(\mathcal{X}, d, \mathcal{P})$ is called a
\textit{vector ultrametric space}. If $\mathcal{P}$ is unital and
normal, then $(\mathcal{X}, d, \mathcal{P})$ is called a
unital-normal vector ultrametric space. }\end{definition}

For any unital-normal vector ultrametric space $(\mathcal{X},
d,\mathcal{P})$ with a convex modular, since
$$d(x,y)\preceq \rho(d(x,y))e \,\,\,\, \mbox{and} \,\,\,\,\, d(y,z)\preceq \rho(d(y,z))e, $$
from (CUM3) we have $$d(x,z)\preceq \max \{\rho(d(x,y)),
\rho(d(y,z)) \}e,$$ and therefore \begin{equation}\label{ineq}
\rho(d(x,z))\leq \max \{\rho(d(x,y)), \rho(d(y,z)) \}.
\end{equation}

For a unital-normal vector ultrametric space $(\mathcal{X},d,
\mathcal{P})$, if $x\in \mathcal{X}$ and $p\in
\mathcal{P}\setminus \{0\}$, the subset
$$B(x;p):=\{y\in \mathcal{X}:\, \rho(d(x,y)) \leq  \rho(p) \},$$ is said to be
     a ball centered at $x$ with radius $p$. Every point of a ball is its center and intersecting balls with comparable radii are comparable with respect to inclusion.  The
unital-normal vector ultrametric space $(\mathcal{X},d,
\mathcal{P})$ is called \textit{spherically complete} if every
chain of balls (with respect to inclusion) has a nonempty
intersection.

\begin{ex}\label{aval}
{\rm Consider the full matrix algebra $\mathbb{M}_n$ over complex
numbers and choose a nonzero positive definite matrix $p$ of
positive cone $\mathcal{P}$ consisting of all positive definite
matrices. \begin{enumerate}
\item For any nonempty set $\mathcal{X}$, define the mapping $d$ by  $$d(x,y)=\left \{\begin{array}{cc}
                                                                                   p & x\neq y \\
                                                                                   0 & x=y.
                                                                                 \end{array}
                                                                                 \right.$$
Then, $d$ is a vector ultrametric on $\mathcal{X}$.
\item Let $(\mathcal{N},\|\cdot \|)$ be a normed space, $(\alpha_n)$ a sequence of positive real numbers decreasing to zero, and
$$\mathcal{X}:=\{x=(x_n)_{n=1}^\infty \in \mathcal{N}: \limsup_{n \rightarrow \infty}\|x_n\|^{\alpha_n}<\infty \}.$$
Now, the mapping $d$ defined by $$d(x,y)=\left \{\begin{array}{cc}
                                                                                   p \limsup_{n\rightarrow \infty} \|x_n-y_n\|^{\alpha_n} & x\neq y \\
                                                                                   0\quad\quad\quad\quad\quad\quad\quad\quad\quad\quad\quad & x=y,
                                                                                 \end{array}
                                                                                 \right.$$
is a vector ultrametric on $\mathcal{X}$.
\item \label{star} Let $\mathcal{A}$ be a C$^\ast$-algebra with positive cone
$\mathcal{P}$ (consisting of the set of all self-adjoint elements
with non-negative spectral values). If $(\mathcal{X},d)$ is an
ultrametric space in the usual sense and $p\in
\mathcal{P}\setminus \{0\}$, then the mapping
$$(x,y)\rightarrow d(x,y)p
\,\,\,\,\,\,\,\,\,\,\,\,\,\,\,\,\,\,\,\, (x,y\in \mathcal{X}),$$
is a vector ultrametric on $\mathcal{X}$.
\end{enumerate}}

\end{ex}

 The next  example generalizes the idea given in the previous example.
\begin{ex} {\rm
Let $\mathcal{A}_\rho$ be a complete modular space with the  cone
$\mathcal{P}$. For usual ultra metric space $(\mathcal{X},d)$ and
$p\in \mathcal{P}\setminus \{0\}$, the mapping
$$(x,y)\rightarrow d(x,y)p
\,\,\,\,\,\,\,\,\,\,\,\,\,\,\,\,\,\,\,\, (x,y\in \mathcal{X}),$$
is a vector ultrametric on $\mathcal{X}$.}
\end{ex}

It is clear that the cones given in Example \ref{aval} are normal
and the cone in \ref{star} of the same example  is also unital
(see, e.g., \cite{conway})

\begin{ex}{\rm  Consider the Euclidean  space  $\mathbb{R}^2$ with the lexicographical order $\preceq$ (i.e., $(a,b)\preceq (a',b'$) if
$a < a'$ or [$a = a'$ and $b\leq b'$]) . Then, it is clear that
$\mathcal{P}=\{x\in \mathbb{R}^2: x\succeq 0\}$ is not normal. For
any nonempty set $\mathcal{X}$ equipped with the mapping
$$d(x,y)=\left \{\begin{array}{cc}
                                                                                   u & x\neq y \\
                                                                                   0 & x=y,
                                                                                 \end{array}
                                                                                 \right.$$
                                                                                 where $u\in \mathcal{P}$  is a fixed element,
                                                                                 we obtain a non-normal and unital vector ultrametric space. In
                                                                                 fact, $(a,b)\preceq \|(a,b)\|(1,1)$, for every $(a,b)\in \mathbb{R}^2$.}
\end{ex}


\section{Main Theorem}
We recall that a correspondence  $\varphi$  on a set $\Omega$,
denoted by $\varphi: \Omega \twoheadrightarrow \Omega$,  assigns
to each $w$ in $\Omega$ a (nonempty) subset $\varphi(w)$ of
$\Omega$. For any subset $C$
 of $\Omega$ and correspondence  $\varphi:C\twoheadrightarrow \Omega$, an
element $w\in C$ is said to be a fixed point if $w\in \varphi(w)$.

By a convergent sequence $(x_n)_{n=1}^\infty$ in vector
ultrametric space $(\mathcal{X}, d, \mathcal{P})$, we mean that
there exists an element $x\in \mathcal{X}$ such that
$\rho(d(x_n,x))\rightarrow 0$ as $n\rightarrow \infty$. It is not
difficult to see that for any unital-normal vector ultrametric
space $(\mathcal{X}, d,\mathcal{P})$, the vector ultrametric $d$
is jointly continuous, i.e, if $x_n\rightarrow x$ and
$y_n\rightarrow y$, then $d(x_n,y_n)\rightarrow d(x,y)$.

\noindent We also say that a subset $G$ of $(\mathcal{X},d, \mathcal{P})$ is compact if every sequence in $G$ has a convergent subsequence in $G$. In the following by $\varphi:\mathcal{X}\twoheadrightarrow c(\mathcal{X})$ we mean that $\varphi$ is a correspondence  with compact values.\\

\noindent\textbf{Theorem} $\,$ Let $(\mathcal{X},d, \mathcal{P})$
be a spherically complete unital-normal vector ultrametric space
and $\varphi:\mathcal{X}\twoheadrightarrow c(\mathcal{X})$. If
for every $x,y\in \mathcal{X}$, $x\neq y$,  and $p\in \varphi(x)$
there exists $q\in \varphi(y)$ such that
\begin{equation}\label{assump} \rho( d(p,q)) < \max \{\rho( d(x,p)),\rho( d(x,y)), \rho( d(y,q)) \},
\end{equation}
 then there exists $g\in \mathcal{X}$ such that
$g\in \varphi(g)$.\\

\noindent \emph{Proof.} Let  $$\Gamma=\{B_{(a,p)}\mid a\in
\mathcal{X},\, p\in \varphi(a)\},$$ where $B_{(a,p)}=B(a;d(a,p))$.
Consider the partial order $\,\sqsubseteq \,$on $\Gamma$ defined
by
$$B_{(a,p)}\sqsubseteq B_{(b,q)}\,\,\,\,\,\,\, \mbox{iff} \,\,\,\,\,\,\,\,
B_{(b,q)}\subseteq B_{(a,p)},$$ where $a,b\in \mathcal{X}$, $p\in
\varphi(a)$, and $q\in \varphi(b)$. If $\Gamma'$ is any chain in
$\Gamma$, then the spherically completeness of $\mathcal{X}$
implies that the intersection $\Omega$ of elements of $\Gamma'$
is nonempty. Choose $c\in \Omega$ and $B_{(a,p)} \in \Gamma'$. If
$x\in B_{(c,q)}$, where $q\in \varphi(c)$ and satisfies
(\ref{assump}) then
$$\rho( d(x,c))\leq \rho( d(c,q)) \leq \max \{\rho( d(c,a)),\rho( d(a,p)),\rho( d(p,q))\},$$
and since $\rho( d(c,a))\leq\rho( d(a,p))$ (because of $c\in
B_{(a,p)}$), we get
\begin{equation}\label{max} \rho( d(x,c) )
\leq \max \{\rho( d(a,p)) ,\rho( d(p,q))\}.
\end{equation}
 We claim that
$\rho( d(x,c)) \leq \rho( d(a,p)).$
 If $ \rho( d(p,q)) \leq \rho( d(a,p))$, then the inequality is clear. If, otherwise $ \rho( d(p,q))
> \rho( d(a,p))$, then from (\ref{max}) we obtain
$$\rho( d(x,c)) \leq \rho( d(p,q)).$$
From (\ref{assump}) it follows that
$$\rho( d(x,c)) < \max \{\rho( d(a,p)),\rho( d(a,c)), \rho( d(c,q)) \},$$
and hence
$$\rho( d(x,c)) < \max \{\rho( d(a,p)), \rho( d(c,q)) \}.$$
Now, if   $\rho( d(a,p)) < \rho( d(c,q))$, then
$$\rho( d(c,q)) \leq  \max \{\rho( d(c,a),\rho( d(a,p)), \rho( d(p,q)) \},$$
that is,  $$\rho( d(c,q)) \leq   \rho( d(p,q)) ,$$ and so from
(\ref{assump}) we get the contradiction
$\rho(d(p,q))<\rho(d(p,q))$. Therefore $$\rho( d(x,c)) \leq \rho(
d(a,p)),$$ and because $B_{(a,p)}=B(c;d(a,p))$, it implies that
$$\rho( d(x,a)) \leq \rho( d(a,p)).$$ That is, $x\in B_{(a,p)}$,
and consequently $B_{(c,q)}\subseteq B_{(a,p)}$. Now,
$$ \inf_{q\in \varphi(c)}\rho(d(c,q))=\rho(d(c,\tilde{q})),$$ for some $\tilde{q}\in \varphi(c)$ (because of (\ref{ineq}) and $\Delta_2$-condition). If $\rho(d(c,\tilde{q}))=0$, then $c\in \varphi(c)$. Otherwise, $B_{(c,\tilde{q})}$ is
 an upper bound for the chain $\Gamma'$. Therefore, by Zorn's
 lemma $\Gamma$ admits a maximal element $B_{(g,w)}$, where $g\in
 \mathcal{X}$ and $w\in \varphi(g)$. We show that $g\in \varphi(g)$. Suppose
 on the contrary that $g\notin \varphi(g)$. Then, by (\ref{assump}), setting $x=g$ and $y=p=w\in \varphi(g)$,
 there exists $s\in \varphi(w)$ such that
 $$\rho(d(s,w))< \max \{\rho( d(g,w)),  \rho(d(w,s))\}$$
 and therefore
 \begin{equation}\label{cont}
 \rho(d(s,w))< \rho( d(g,w)).
 \end{equation}
 On the other hand, from the maximality of $B_{(g,w)}$ and that $ w\in B_{(g,w)}$, we have
$$ B_{(g,w)}\subseteq B_{(w,s)}= B(g;d(w,s)),$$
and so
 $$ \rho(d(w,g)\leq \rho(d(w,s)),$$ which contradicts (\ref{cont}).

The following  corollaries obtain  immediately from preceding
theorem. We suppose that $(\mathcal{X},d, \mathcal{P})$, $\gamma$,
and $\varphi$ are as given in the previous  theorem.

\begin{corollary}
If for  every $x,y\in \mathcal{X}$, $x\neq y$,  and $p\in
\varphi(x)$ there exists  $q\in \varphi(y)$ such that
$$\rho( d(p,q)) < \rho( d(x,y)),$$
then there exists $g\in \mathcal{X}$ such that $g\in \varphi(g)$.
\end{corollary}

\begin{corollary}
If for every $x,y\in \mathcal{X}$, $x\neq y$,  and $p\in
\varphi(x)$ there exists $q\in \varphi(y)$ such that
$$\rho( d(p,q)) < \max \{\rho( d(x,p)),\rho( d(x,y)), \rho( d(y,q)) \},$$ then $\varphi$ has a fixed point.
\end{corollary}

\begin{corollary}
If for every $x,y\in \mathcal{X}$, $x\neq y$,  and $p\in
\varphi(x)$ there exists $q\in \varphi(y)$ such that
$$\rho( d(p,q)) < \rho( d(x,y)),$$ then $\varphi$ has a fixed point.
\end{corollary}

As seen, the last corollary generalizes Theorem 1 in \cite{petalas}.\\

\noindent \textbf{Acknowledgment} The author would like to thank
professor C.G. Petalas for the valuable comments on this note.


\end{document}